\documentclass{amsart}

\usepackage{amssymb}
\newtheorem{thm}{Theorem}[section]
\newtheorem{cor}[thm]{Corollary}
\newtheorem{lem}[thm]{Lemma}
\newtheorem{prop}[thm]{Proposition}
\theoremstyle{definition}
\newtheorem{defn}[thm]{Definition}
\theoremstyle{remark}
\newtheorem{rem}[thm]{Remark}
\newtheorem{ex}[thm]{Example}
\numberwithin{equation}{section}

\def\clap#1{\hbox to 0pt{\hss#1\hss}}

        \def\mathllap{\mathpalette\mathllapinternal}
        \def\mathrlap{\mathpalette\mathrlapinternal}
        
        \def\mathllapinternal#1#2{%
        \llap{$\mathsurround=0pt#1{#2}$}}
        \def\mathrlapinternal#1#2{%
        \rlap{$\mathsurround=0pt#1{#2}$}}

\begin{document}


\title{Non-commutative duplicates of finite sets}%
\author{Claude CIBILS}%
\address{Universit\'{e} Montpellier 2\\Institut de Math\'{e}matiques et de Mod\'{e}lisation de Montpellier\\
F-34095 Montpellier Cedex 5\\France}%
\email{Claude.Cibils@math.univ-montp2.fr}%



\begin{abstract}
We classify the twisted tensor products of a finite set algebra
with a two elements set algebra using colored quivers obtained
through considerations analogous to Ore extensions. This provides
also a classification of entwining structures between a finite set
algebra and the grouplike coalgebra on two elements. The resulting
$2$-nilpotent algebras have particular features with respect to
Hochschild (co)homology and cyclic homology.
\end{abstract} \maketitle


\section{Introduction}
A geometric category can be studied through the category of
commutative algebras of functions on the spaces with values in a a
commutative ring or a field, where the functions have properties
according to the geometrical nature of the spaces; both categories
are usually equivalent through this contravariant functor. Using
this equivalence the algebraic approach is a substitute for the
geometric one, keeping the geometric origine as a background. In
this context the tensor product of algebras corresponds to the
product of the spaces. The non-commutative geometry approach uses
non-commutative algebras as substitutes for the spaces, providing
a setting corresponding to virtual non-commutative spaces.

The set algebra of $k$ valued functions over a a finite set $E$ is
a commutative semi-simple basic finite dimensional $k$-algebra,
isomorphic to $k\times\mathrm{}\cdots\times\mathrm{}k$ where the
number of copies of $k$ equals the cardinality of $E$. The tensor
product of algebras $k^{E}\!\otimes_k k^{F}$ corresponds to the
cartesian product $E\times F$. The non-commutative counterpart
consist on twisted tensor products of $k^{E}$ and $k^{F}$
corresponding to non-commutative {virtual} cartesian products of
the sets $E$ and $F$. Recall that a twisted tensor product of two
algebras is an algebra structure on the tensor product of the
underlying vector spaces such that the natural inclusions of the
original algebras are algebra maps.  This structure has been
considered by P. Cartier  when lecturing in
 Paris, \cite{car}.
Several authors have also considered such product in different
contexts, see the work of
 J. Baez \cite{ba},  A. Cap, H. Schichl, J. Van\v zura \cite{ca} and G. Maltsiniotis \cite{ma}.
   Note that twisted
tensor products between an algebra and the polynomial algebra in
one variable corresponds precisely to Ore extensions, see for
instance \cite {ka} for the definition and use of Ore extensions.

The problem of classifying structures arising this way is
difficult. Quite diverse non semi-simple one parameter families of
algebras can occur. In this paper we use the equivalence between
finite sets and semi-simple basic algebras over a field $k$ in
order to classify the non-commutative duplicates of a set. More
precisely the usual duplicate of a set $E$ is the disjoint union
of two copies of $E$, namely $E \times\mathrm{}\{a,b\}$ which has
algebra set
$$k^{E \times \{a,b\}} \;= \;k^{E} \otimes \;k^{\{a,b\}} = \;k^{E}
\oplus \;k^{E}.$$

Our purpose in this paper is to concentrate in this case
corresponding to non-commutative duplicates of a finite set $E$.

We recall first some facts concerning twisted tensor products of
algebras and how they are obtained. Then we provide a
classification of the twisted tensor products of $k^{E}$ and
$k^{\{a,b\}}$ through the set of interlacings which provides an
analogous setting to Ore extensions. The results shows that
non-commutative duplicates of a finite set are provided
 by one-valued quivers with set of vertices $E$ which
has even length oriented cycles if any and provided with a
coloration. The corresponding non semi-simple non-commutative
2-nilpotent algebras obtained through the path algebra
construction produces the set of non-commutative duplicates of
$E$.

This result also provides a classification of entwinings maps
between the algebra $k^E$ and the two-dimensional grouplike
coalgebra $k\{a,b\}$. Indeed, an entwining consists of an algebra,
a coalgebra, and a map satisfying certain properties, see
\cite{br} or \cite[Sec. 2.1]{camizh}. In case the coalgebra is
finitely generated and projective over the ground ring, there is a
bijection between entwinings and interlacings between the algebra
and the dual of the coalgebra, see \cite[Theorem 8]{camizh}.
Finalley notice that $k^{\{a,b\}}$ is the dual of the grouplike
coalgebra $k\{a,b\}$.

 In the last section we study the Hochschild cohomology of
the non-commutative duplicates of a finite set $E$. The main
result is that the Hochschild cohomology is finite dimensional if
and only if the quiver associated to the situation has no oriented
cycles other than loops, regardless of the coloration. In that
case the only possible non zero vector spaces lies in degrees $0$
and $1$.

We also compute cyclic homology of this algebras in characteristic
zero showing that the dimension is constant in even degrees. In
odd degrees, the dimension of the cyclic homology counts the
number of connected components having a proper circuit of length
dividing the degree augmented by one one.

Concerning Hochschild homology, we obtain that its dimension is
finite if and only if the quiver has no oriented cycles, in which
case all the homology vector spaces vanishes except in degree $0$
where the dimension equals the dimension of the even cyclic
homology.

\section{Twisted tensor products and interlacings}

Let $k$ be a field. We consider $k$-algebras which are associative
and unitary algebras over $k$. We recall some basic facts
concerning twisted products of algebras.

Let $A$ and $B$ be $k$-algebras and let $A\otimes_k B$ be the
tensor product of the underlying $k$-vector spaces. A twisted
tensor product algebra of $A$ and $B$ is a $k$-algebra structure
on $A\otimes_k B$ such that the canonical inclusions $a\mapsto
a\otimes 1$ and $b \mapsto 1\otimes b$ are $k$-algebra maps. Note
that we do not assume any extra structure or possible actions on
the original algebras which could provide the algebra structure of
the twisted tensor product. Of course an action or a grading of a
group on an algebra provide examples of twisted products, see
\ref{actions}.

\begin{lem}
Let $\Lambda = A\otimes_k B$ be a twisted tensor product algebra.
Then $$(a\otimes 1) (1\otimes b) = a \otimes b.$$
\end{lem}

 {\bf Proof}
    The product map $m : \Lambda\otimes_k
    \Lambda  \rightarrow \Lambda$ is a
    $\Lambda$-bimodule map and the restriction  to $\bigl( (A\otimes 1)\otimes(1\otimes
    B)\bigr)$ is an $A - B$ - bimodule map. Moreover $\bigl( (A\otimes 1)\otimes(1\otimes
    B)\bigr)$ is a free $A - B$ bimodule generated by $\bigl((1\otimes 1)\otimes(1\otimes 1)\bigr)$ and
we have    $m\bigl ((1\otimes 1)\otimes(1\otimes 1)\bigr) =
\bigl(1\otimes
    1).$ Consequently
 \begin{eqnarray*}
   m\;\bigl( (a\otimes 1)\otimes(1\otimes
    b)\bigr)    &=& m\;\bigl( a\; (1\otimes 1)\otimes(1\otimes
    1)\; b\bigr)\\
&=& a\;\bigl( (1\otimes 1)\otimes(1\otimes
    1)\bigr)\;b\\
&=& a \;m\;\bigl( 1\otimes 1)\;b\\
&=&a \otimes \;b.
\end{eqnarray*}

As a consequence the structure of a twisted product is determined
by the value of the products of elements of $B$ with elements of
$A$. Let $\tau : B \otimes A\to A \otimes B$ be the map given by
$\tau (b \otimes a) = (1 \otimes b) (a \otimes 1)$, in other words
$\tau$ provides the product of $b$ and $a$ inside the twisted
tensor product of algebras.

 \begin{prop}
   The map $\tau$ verifies $\tau ( b \otimes 1) = 1 \otimes b$,
    $\tau (1 \otimes a) = a \otimes 1$ and the  braiding
    conditions, namely the diagrams below commute (tensor signs are omitted or replaced by commas)
\end{prop}

 $$ \begin{array}{ccccc}
BBA&\stackrel{\scriptstyle1_B,\tau}{\longrightarrow}&BAB&
\stackrel{\scriptstyle\tau,1_B}{\longrightarrow}&ABB\\
\mathllap{\scriptstyle m_{B},1_A}\left\downarrow\vbox to
1.5em{}\right.&&
&&\left\downarrow\vbox to 1.5em{}\right.\mathrlap{\scriptstyle 1_A,m_{B}}\\
BA&\multicolumn{3}{c}{\stackrel{\tau}{\hbox to
8em{\rightarrowfill}}}&AB

  \end{array}$$

\quad
$$
  \begin{array}{ccccc}
BAA&\stackrel{\scriptstyle \tau,1_A}{\longrightarrow}&ABA&
\stackrel{\scriptstyle1_A,\tau}{\longrightarrow}&AAB\\
\mathllap{\scriptstyle 1_B, m_A}\left\downarrow\vbox to
1.5em{}\right.&&
&&\left\downarrow\vbox to 1.5em{}\right.\mathrlap{\scriptstyle m_A,1_{B}}\\
BA&\multicolumn{3}{c}{\stackrel{\tau}{\hbox to
8em{\rightarrowfill}}}&AB
  \end{array}
$$

 {\bf Proof:} The first diagram correspond to the associativity
condition of the product for triples $(b,b',a)$ namely
$(b\otimes1) (b' \otimes 1) ( a \otimes 1)$ while the second
corresponds to triples $(b,a,a').$

\begin{defn}
An interlacing is a linear map map $\tau : BA \longrightarrow AB$
verifying the preceeding properties.
\end{defn}

Note that an interlacing is also called a factorization structure,
see for instance \cite{maj}. The following result is now obvious :

\begin{thm}
 Let $A$ and $B$ be $k$-algebras. The twisted tensor
 product structures on $A\otimes_k B$ are in one to one
 correspondence with interlacings.
\end{thm}

 \begin{ex}
 The usual tensor product of algebras is obtained through the
 trivial interlacing, namely the flip map $\tau ( b \otimes a ) =
 a \otimes b.$
\end{ex}

\begin{ex}\label{actions}
Let $G$ be a group acting by automorphisms on a $k$ algebra $A$
and let $\tau$ be the map given by $\tau (s \otimes a ) = s (a)
\otimes s$. Then $\tau$ is an interlacing which provides the well
known skew group algebra $A [G]$. As usual this can be
generalized, replacing $kG$ by a Hopf algebra $H$ and the action
by an $H$-module algebra structure on $B$.
\end{ex}

\medskip

Our purpose is to investigate the set of twisted tensor products
between two algebras in order to provide a classification and
properties relatively to the original algebras. As quoted in the
Introduction, in this paper we will consider the case where $B =
k^{\{a,b\}}$ is the set algebra of two elements $\{a,b\}$. Note
that the trivial interlacing provides the usual tensor product $A
\otimes k^{\{a,b\}} = A\oplus A$. Note also that the two elements
algebra $k^{\{a,b\}}$ is isomorphic to the algebra of truncated
polynomials $k[X] / (X^{2} - X)$ where $a$ correspond to $X$ and
$b$ to $1-X.$

In order to classify interlacings we consider the following set.

\begin{defn}
Let $Y_{A}$ be the set of couples $(f,\delta)$ such that $f$ is an
endomorphism of $A$ and
$A\stackrel{\delta}{\longrightarrow}{}^{f}\!\!A$ is an idempotent
derivation verifying $$f = f^{2} + \delta f + f \delta.$$
\end{defn}

 \begin{rem}
 The notation ${}^{f}\!\!A$ stands for the A-bimodule given by the
 vector space A with left action modified by $f$, namely $a \cdot a' =
 f(a)a'.$ The right action is the standard one given by the
 product of $A.$
\end{rem}

\begin{rem}
The data above can be interpreted as a sort of Ore extensions
associated to the quotient of the polynomial algebra in one
variable that we consider. One can prove that Ore extensions (see
for instance \cite{ka}) are precisely twisted tensor products of a
$k$-algebra with the polynomial algebra in one variable $k[x]$.
The proof of this fact can be performed following the lines of the
next Proposition.
\end{rem}

 \begin{prop}
The set $Y_{A}$  is in one-to-one correspondence with the set  of
interlacings between the two elements set algebra $k[X]/(X^{2} -
X)$ and a $k-$algebra  $A$.
\end{prop}

 {\bf Proof:} Consider the $k$-module $A \otimes [k[X]/(X^{2}-X)] =
A[X]/(X^{2}- X).$ An interlacing $\tau$ is determined by the
values $\tau (X, a)$ corresponding to the product $X a.$ We put
$\tau (X, a) = X a = \delta(a) + f(a) X.$ The braiding conditions
provides the equalities $X (X a) = X a $  and $X(ab) =
    (X a) b.$
The first one translates into $\delta^{2} = \delta$ and $f =
    f^{2} + \delta f + f \delta$ while the second provides
     that $f$ is an endomorphism of $A$ and $\delta$ is a
derivation with coefficients in ${}^f\!A.$ Conversely it is clear
that those conditions insures that the
    map $\tau$ defined by the above formulae is an interlacing.
Note that in this setting the trivial interlacing corresponds to
the identity endomorphism
 and the zero derivation.

\begin{rem}
One could expect a cohomological interpretation of $Y_{a}$. It
appears that $Y_A$ do not seem to have any additional structure,
in particular no natural composition law of elements of $Y_A$ seem
to exist. In the next section it will become clear that interior
 derivations do not provide in general trivial tensor twisted products.
\end{rem}

\begin{defn}
The set of \emph{$2$-interlacings} of an algebra $A$ is the set of
interlacings between $A$ and $k[X]/(X^{2} - X)$. We have proved
that this set is in one-to-one correspondence with $Y_{A}$.
\end{defn}

\section{Set algebras and $2$-interlacings}

Our first purpose in this section is to describe $Y_{{k}^{E}}$
where $k^{E}$ is the set algebra of a finite set $E$. The second
purpose is to classify the family of twisted tensor algebras
obtained through the $2$-interlacings. Recall that $k^{E}$ is the
$k$-algebra of functions on $E$, namely $k^{E} = \{ a : E
\longrightarrow k \}$ where the product is given by $(aa')(x) = a
(x) a' (x).$ Of course $k^{E}$ is isomorphic to the vector space
with basis the set $E$ and componentwise product

$$(\sum_{x\in E} a_{x}x)(\sum_{x\in E} a'_{x}x)
= \sum_{x\in E} (a_{x}a'_{x})x.$$

In other words $k^{E}$ is a product of copies of $k$ indexed by
$E$, and $E$ is a complete set of primitive orthogonal idempotents
of $k^{E}.$

It is well known that $A = k^{E}$ is cohomologicaly trivial, which
means that the Hochschild cohomology of $A$ with coefficients in
any $A$-bimodule vanish in positive degrees. In particular any
derivation is interior, this fact will provide an exhaustive way
for describing the set of  $2$-interlacings. For the convenience
of the reader as well as in order to introduce notation we provide
a proof of this result.

Let $A$ be a $k$-algebra and M an $A$-bimodule. A
\emph{derivation} is a $k$-linear map  $\delta: A\to M$ such that
$\delta(a a') = a \delta (a') + \delta\, (a) a'.$ An
\emph{interior} derivation is associated to each element $m$ of
$M$ by the formula $\delta_{m} (a) = am - ma.$ One can check that
interior derivations are indeed derivations.

The following result is well known and easy to prove, it provides
the classification of $A$-bimodules when $A$ is a finite set
algebra over a field.

 \begin{prop}
 Let $k$ be a field, $E$ be a set and $A = k^{E}$ be the set
 algebra. Any finitely generated $A$-bimodule is
  isomorphic to a direct sum of simple modules. The complete list
 of simple modules up to isomorphism is $\{_{v}k_{u} \}_{u, v\in E}$
 where $_{v}k_{u}$ is a one dimensional vector space with identity
 actions of $v$ on the left, of $u$ on the right, and zero actions of
 others set elements.

\end{prop}

According to this Proposition it suffices to consider simple
bimodule of coefficients in order to prove that any derivation is
inner.

 \begin{prop}
 Let $A = k^{E}$ be a set algebra, let $_{v}k_{u}$ be a simple
 bimodule and let $\delta : A \longrightarrow {}_{v}k_{u}$ be a
 derivation. If $u \neq v$ the derivation is inner. In case $u = v$
 we have $\delta = 0$.
\end{prop}

 {\bf Proof:} We prove first that the space of derivations is
one-dimensional if $u \neq v$, and is zero if $u = v.$ Let $e\in
{}_{v}k_{u}$ be a fixed non zero element and for each $x\in E$ let
$\lambda_x\in k$ defined by $\delta(x) =
 \lambda_x e$. We assert that $\lambda_x = 0$ if
 $x\neq u$ or $x \neq v.$ Indeed,
 $$\lambda_{x}e = \delta (x) = \delta (x^{2}) = x
 \delta (x) + \delta (x) x = \lambda_x {x e} + \lambda_x e x =
 0+0=0.$$
If $u = v$ we have $\lambda_{u} e = \delta (u) = 2 \lambda_{u} e$
which implies $\lambda_{u} = 0.$ If $u \neq v$ then  $$0 = \delta
(vu) = v \delta (u) + \delta (v) u =
 \lambda_{u}ve + \lambda_{v}e u
 = (\lambda_{u} + \lambda _{v}) e.$$
 We consider now the interior
derivation $\delta_{e}$ given by $\delta_e (x) = xe - ex.$ Of
course we have  $\delta_{e} (x)=0$ if $x \neq u$ or $x \neq v.$
Assuming $u \neq v$, we obtain $\delta_{e} (v) = e$ while
$\delta_{e} (u) = - e.$ Consequently the space of interior
derivations is also one dimensional and every derivation is inner.

Towards the description of the set of $2$-interlacings, we recall
that each algebra endomorphism of $k^{E}$ is determined by a set
map $\varphi : E \longrightarrow E.$ Actually this is a special
case of the anti-equivalence between the category of finite sets
and the category of semi-simple basic and commutative algebras.
More precisely let $f$ be an algebra endomorphism of $k^{E}$.
There exist a unique set map $\varphi : E \longrightarrow E$ such
that for each $e\in E$ we have

$$f (e) = \mathop \Sigma _{\{x\mid\varphi(x) = e\}}x.$$

\begin{lem}
Let $f$ an algebra endomorphism of $A = k^{E}$ given by a set map
$\varphi : E \longrightarrow E.$ Let $\delta : A \longrightarrow
{}^{f}\!\!A$ be a derivation. There exist $a\in A$ such that for
each  $e\in E$ we have $
  \delta(e) = (f(e) - e) a
   = \Sigma_{\varphi(x)=e}  a_{x} x - a_{e}e.
$ Moreover $$\delta^{2}(e) = \sum_{\varphi(y) = e} a_{y} a
_{\varphi (y)}y \ -
   \sum_{\varphi(x) = e} a_{x}(a_{x} + a_{e}) x \ + \ a^{2}_{e} e.$$
\end{lem}

\begin{rem}
The element $a$ of the Lemma is uniquely determined once
normalized  at loop elements of $E$, namely $a_{e} = 0$ if $f(e) =
e.$ We call if the determining element of $\delta$.
\end{rem}

 {\bf Proof:} Since $\delta$ is interior there exist $a$ such that
   $\delta(e)=e \cdot a - a_{e}e
 = f(e)a - ae\\
  = (f(e) - e)a.$
The fact that $a$ is unique follows at once from the previous
considerations. Note also that if $u\neq v$ we have $H^0(A,
{}_{v}k_{u}) =0$.

 Now
\begin{eqnarray*}
\delta^{2}(e) &=& \sum_{\varphi(x) = e} a_{x}\delta(x) - a_{e}\delta(e)\\
   &=& \sum_{\varphi(x) = e} a_{x} \left(\sum_{\varphi(y) = x} a_{y}y - a_{x}x\right) -
   a_{e} \left(\sum_{\varphi(x) = e} a_{x}x - a_{e}e\right)\\
   &=&  \sum_{\varphi(y) = e} a_{y} a _{\varphi (y)}y \ -
   \sum_{\varphi(x) = e} a_{x}(a_{x} + a_{e}) x \  + \ a^{2}_{e} e.\\
\end{eqnarray*}

In order to describe the idempotent derivations  $\delta : A
\longrightarrow {}^{f}\!\!A$  it is useful to introduce the quiver
of the endomorphism $f$ given by a set map $\varphi : E
\longrightarrow E.$ Recall that a quiver $Q$ is a finite oriented
graph with set of vertices $Q_{0}$,  set of arrows $Q_{1}$ and two
maps $s$, $t$ : $Q_1 \longrightarrow Q_{0}$ providing each arrow
with a source and a target vertex.

\begin{defn}
Let $f$ be an endomorphism of the set algebra $A = k^{E}$ given by
a set map $\varphi.$ The quiver $Q_{f}$ of $f$ has set of vertices
$E$ and an arrow from $x$ to $\varphi(x)$ for each $x\in E.$ Two
arrows  $b$ and $a$ are concatenated if $s(b)=t(a)$
\end{defn}

\begin{rem}
Quivers obtained this way are precisely \emph{one-valued quivers},
that is each vertex of the quiver is the source of exactly one
arrow, accordingly to the definition of a set map. \end{rem}

\begin{defn}
An \emph{oriented cycle} of a quiver is a sequence of concatenated
arrow such that the source of
 the first arrow coincides with the target of the last one.
 Its \emph{length} is the number of arrows involved. A \emph{loop} is
 an oriented cycle of length one, namely an arrow with the same source and
 target vertices. An oriented cycle is \emph{proper} if it is not the iteration of a
 strictly smaller length  oriented cycle. A \emph{loop vertex} is
 a vertex where a loop has its source and target vertices.
\end{defn}

Each connected component of a one-valued quiver has precisely one
proper oriented cycle, which can be a loop.

 Let
$R$ be a connected component of the quiver $Q_{f}.$ Let $\delta :
k^{E}\longrightarrow {}^{f}(k^{E})$ be a derivation with
determining element $a\in k^{E}.$ Our first purpose is to describe
those $a$ such that $\delta^{2} = \delta.$

\begin{lem}
Let $A = k^{E}$ be a finite set algebra, $f$ an endomorphism of
$A$ with set map $\varphi$ and quiver $Q_{f} $ and let $A
\longrightarrow {}^{f}\!\!A$ be an idempotent derivation with
determining element $a\in k^{E}$. Let $u\longrightarrow v$ be an
arrow of $Q_{f}$ with $u \neq v$ where $v$ is a non loop vertex
and such that there is no arrow backwards. Note that $u$ cannot be
a loop vertex since $Q_{f}$ is one-valued. Then $a_{u}, a_{v}
\in\{{-1, 0}\}$ and $a_{u} a_{v} = 0.$
\end{lem}

 {\bf Proof:} At the vertex $v$ we have the formula
$$0 = (\delta^{2} - \delta)
(v) =\sum_{\varphi^{2}_{y} = v}  a_y a_{\varphi(y)} y \ -
\sum_{\varphi(x) = v} a_{x} (a_{x} + a_{v} + 1)x \ + \ a_{v}(a_{v}
+ 1)v.$$ The $v$-coefficient of this sum is $a_{v}(a_{v} + 1)$,
then $a_{v}\in \{0, -1\}.$ The $u$-coefficient is $-a_{u}(a_{u} +
a_{v} + 1).$ If $a_{v} = -1$ then $a_{u} = 0.$ If $a_{v} = 0$ then
$a_{u}\in\{0, -1\}.$

Next we define pre-colorations of $Q_f$, which will correspond to
idempotent derivations.

\begin{defn}
Let $Q_{f}$ be the quiver of an algebra endomorphism $f$ provided
by a set map $\varphi : E \longrightarrow E.$ A
\emph{pre-coloration} of its vertices is an element $a$ of $k^{E}$
verifying the following conditions:

\begin{enumerate}
    \item In case of a connected component of $Q_{f}$ reduced to a round trip
quiver $u\cdot\!\rightleftarrows\! \cdot v$ the colors $a_{u}$ and
$a_{v}$ verify either $a_{u} + a_{v} + 1 = 0$ or $a_{u} = a _{v} =
0.$
    \item   For a connected component different from the round trip
 quiver:
\begin{enumerate}
    \item  In case $e$ is a non loop vertex then $a_{e} \in \ \{{0, -1}\}.$
    \item  For each arrow the product of the
    colors at its source and target is $0$.
    \item  For a loop vertex the color is irrelevant, we normalize it at 0.
\end{enumerate}
\end{enumerate}

\end{defn}

 \begin{prop}
Let $A = k^{E}$ be a finite set algebra, $f$ an endomorphism of
$A$ with set map $\varphi$ and quiver $Q_{f} $ and let $A
\longrightarrow {}^{f}\!\!A$ be an idempotent derivation with
determining element $a\in k^{E}$. Then $a$ is a pre-coloration of
$Q_{f}$ on this connected component.
\end{prop}

 {\bf Proof:}
We consider a connected component $R$ of $Q_f$. Assume first there
is a loop vertex $e$ in $R$ and consider the equation
 $(\delta^{2} - \delta)(e) = 0.$
 The coefficient of $e$ is $a^{2}_{e} - a_{e} ( 2 a_{e} + 1 ) +
a^{2}_{e} + a_{e}$ which is $0$ for any value of $a_{e}.$ If the
connected component $R$ is reduced to $e$ we have that $a$ is a
pre-coloration. Otherwise let $x$ such that $\varphi(x) = e$ and
$x \neq e.$ The coefficient of $x$ in the equation $(\delta^{2} -
\delta) (e) = 0$ is
$$a_{x}a_{e} - a_{x} (a_{x} + a_{e} +  1) =
a_{x} (a_{x} + 1).$$
 Then $a_{x}\in\{0, -1\}.$ In case $R$ do not contain the round trip quiver the conclusion
follows from the preceding lemma.

If $R$ contains the round trip quiver having vertices $u$ and $v$
the equation $(\delta^{2} - \delta) (u) = 0$ provides the
$v$-coefficient
$$-a_{v} (a_{v} + a_{u} + 1)$$
and $u$-coefficient
$$a_{u}a_{v} + a_{u}(a_{u} + 1) = a_{v}(a_{v} + a_{u} + 1).$$

It follows that $a_{u} + a_{v} + 1 = 0$ or $a_{u} = a_{v} = 0.$

In case $R$ is not reduced to the round trip quiver there is an
arrow arriving to $u$ or $v$ -- we assume $u$ is this vertex
without lost of generality -- coming from a vertex $w$ which
cannot be a loop vertex since each vertex is the source of exactly
one arrow. For the same reason this arrow cannot be part of a
round trip quiver. The Lemma applies and we infer $a_{w}, a_{u}
\in \{0, -1\}$ and $ a_{u} a_{w} = 0.$ If $a_{u} = 0$ the
equations of the round trip quiver above provides $a_{v} = -1$ or
$a_{v} = 0$. If $a_{u} = -1$ then $a_{v} = 0$. In both cases we
obtain a pre-coloration.

\begin{rem}
From the proof of the result it is clear that conversely a
pre-coloration of $Q_{f}$ provides an idempotent derivation $k^{E}
\longrightarrow  {}^{f}\!(k^{E}).$
\end{rem}

Finally we need to describe the pre-colorations corresponding to
idempotent derivations verifying $f = f^2 + \delta f + f\delta$.

\begin{defn}
A \emph{coloration} of $Q_{f}$ is an element $a = \sum_{x\in E}
a_{x} x \in k_{E}$ such that
 \begin{enumerate}
    \item  For a connected component reduced to the round trip
    quiver
    $u\cdot\!\rightleftarrows\!\cdot v$ the coefficients $a_{u}$ and $a_{v}$
    verify $a_{u} + a_{v} + 1 = 0$.
    \item  For other connected components:
    \begin{enumerate}
        \item  In case $e$ is a non loop vertex, $a_{e} \in
        \{0, -1\}$
        \item  For each arrow having no loop  vertex target, one extremity value is
        $0$ and  the other is $-1.$
        \item  At a loop vertex the value of $a$ is irrelevant,
        we normalise it at $0.$
    \end{enumerate}
\end{enumerate}
\end{defn}

 Note  that a coloration is a pre-coloration.

\begin{rem}\label{description}
Not any  quiver $Q_{f}$
 admits a coloration, clearly the length of non loop proper oriented cycles
 has to be even. For a connected
 component having even length proper oriented cycles, precisely two
 colorations exists. Otherwise a coloration of a connected component with a loop is completely
 determined by providing $0$ and $-1$ values on the loop related vertices, namely on
 source of arrows having a loop target vertex.
 \end{rem}

\begin{thm}
Let $E$ be a finite set, $A = k^{E}$, and let $f$ be an
endomorphism of $A$ corresponding to a set map $\varphi$ with
quiver $Q_{f}.$ Let $a \in k^{E}$ be a normalized element (i.e. it
has $0$ value at loop vertices) ant let $\delta: k^{E}
\longrightarrow {}^{f}\!(k^{E})$ be the inner derivation
determined by $a$. The couple $(f, \delta)$ is   a $2$-interlacing
for $A$ if and only if  $\ a$ is a coloration of $Q_{f}$.
\end{thm}

 {\bf Proof:}
We already know that pre-colorations translates the fact that
$\delta^{2} = \delta.$ The condition $f = f^{2} + \delta f +
f\delta$ provides the following
 for each $e\in E$ :

 $$\sum_{\varphi^{2}(y) = e} (a_y + a_{\varphi(y)} + 1 ) y
 \ - \sum_{\varphi_{x} = e} (a_{e} + a_{x} + 1 ) x = 0.$$

 Let $u \longrightarrow v$ be an arrow between different non loop
 vertices. The $u$-coefficient in the above equation for $v$ provides
 $a_{v} + a_{u} + 1 = 0.$
 Since we have $a_{u}, a_{v} \in \{0, -1\}$ we infer that $a_{u}
 \neq a_{v}.$ In case $v$ is a loop vertex the $u$-coefficient in the equation
for $v$ gives $$(a_{u}+ a_{v} + 1) - (a_{v}+ a_{u} + 1) = 0$$ and
there no restriction is inferred for $a_{v}.$ The $v$-coefficient
for the equation of $v$ provides
 $$(a_{v}+ a_{v} + 1) - (a_{v}+ a_{v} + 1) = 0.$$
Finally assume that $R$ is the round trip quiver $u
 \rightleftarrows v$. The $u$-coefficient for the $v$-equation
 gives
 $$(a_{v}+ a_{u} + 1) = 0$$
 which implies that the option $a_{u} = a_{v} = 0$ of a
 pre-coloration do not provide an interlacing. In contrast the condition
 $a_{v}+ a_{u} + 1 = 0$ is maintained.

\begin{thm}

Let $A = k^{E}$ be the algebra of a finite set $E$. The set
$Y_{A}$ of $2$-interlacings of $A$ is in bijection with one-valued
quivers on the set $E$ provided with a coloration.

\end{thm}

\section{Non-commutative duplicates}

The aim of this section is to classify the twisted tensor products
$k^{E} \otimes k^{\{a, b\}}$ where $E$ is a finite set. We will
obtain a family of non semi-simple, non-commutative algebras, with
square zero Jacobson radical.

From the preceding section it appears that we can assume that the
quiver $Q_{f}$ of an endomorphism $f : k^{E} \longrightarrow
k^{E}$ is connected. More precisely let $(f, \delta)$ be an
endomorphism of $k^{E}$ and let $\delta$ be a derivation
$$\delta : k^{E} \longrightarrow
{}^{f}\!(k^{E})$$ verifying $\delta^{2} = \delta$ and $f = f^{2} +
\delta f + f \delta$. Let $k^{E}\otimes_{(f,\delta)}
\left[{k[X]}/{(X^{2} = X)}\right]$ be the corresponding twisted
tensor product.

Let $E^{1}, ..., E^{n}$ be the partition induced on $E$ by the
connected components of $Q_{f}$ and let $f_{i}$ be the
corresponding endomorphism of $k^{E_{i}}$. Clearly $\delta$
decomposes into $\delta_i:k^{E_{i}}  \longrightarrow {}^{f_{i}}
(k^{E_{i}})$ and we have
$$k^{E} \otimes_{(f,\delta)}{k[X]}/{(X^{2}-X)}
\ =\  \prod_{i=1,\cdots, n}\ \
\left(k^{E_{i}}\otimes_{(f_{i},\delta_{i})}
\left[k[X]/(X^{2}-X)\right]\right).$$

First we consider the trivial case of a one element set $E$, the
map $f$ is the identity, the quiver $Q_{f}$ is a loop and $\delta
= 0$. The only twisted tensor product of $k^{E}$ with $k[X]/(X -
X^{2})$ is the trivial one. There is no non trivial duplicate of a
one element set.

The case where $E$ is an arbitrary finite set but $f$ is the
identity reduces to preceding case, use for instance the connected
components of $Q_{f}$ which are all of them loops, or perform the
direct computation since the only derivation is the zero one.

Recall that for a quiver $Q$ the path algebra $kQ$ has a basis
provided by all the oriented paths -- that is finite sequences of
concatenated arrows of $Q$ -- which multiplies as they concatenate
if this is possible, and have zero product otherwise. By
definition the two sided ideal $\langle Q^{2}\rangle$ has a basis
given by all the paths of length greater or equal $2.$ The
quotient $k$-algebra $kQ/\langle Q^{2}\rangle$ has a basis
provided by vertices and arrows. Its Jacobson radical has a basis
given by the arrows and has zero square. It is well known that any
finite dimensional, basic and split radical square zero algebra is
obtained this way.

We assert that twisted tensor products $k^{E} \otimes
\left[k[X]/(X^{2} - X)\right]$ are members of the preceding family
of algebras. In order to prove this assertion we define a quiver
related to a couple $(f, \delta) \in Y_{k^{E}}$. It will appear to
be essentially the opposite quiver $Q_{f}^{\rm{op}}$ of the set
map $\varphi$ defining $f$, except for loop vertices which becomes
two vertices.
\begin{defn}
Let $E$ be a finite set and $f : k^{E} \longrightarrow k^{E}$ be
an algebra map provided by a set map $\varphi : E \longrightarrow
E$ with  quiver $Q_{f}.$ Let $a$ be a coloration of $Q_{f}$ and
let $\delta$ be the corresponding derivation. The related quiver
$Q_{(f,\delta)}$ is obtained by replacing each connected component
$R$ of $Q_f$ as follows

\begin{itemize} \item
If $R$ has no loops replace it by its opposite $R^{\rm{op}}$ which
has same set of vertices while arrows are reversed.

\item If $R$ has a loop vertex $\ell$  (in which case $R$ do not
contain other proper oriented cycles) remove the loop from
$R^{\rm{op}}$ and replace the vertex $\ell$ by two vertices
$\ell_{0}$ and $\ell_{-1}$. Let $\varepsilon = 0 \mbox{ or } -1$.
Each arrow in $R^{\rm{op}}$ from an $\varepsilon$-polarized vertex
to $\ell$
 is replaced by an arrow from the $\varepsilon$-polarized vertex to
$\ell_{\varepsilon}$. Note that this procedure produces two
connected components.
\end{itemize}

\end{defn}

\begin{thm}
Let $E$ be a finite set and let $f : k^{E}\longrightarrow k^{E}$
be an algebra endomorphism given by a set map $\varphi$ with a
connected quiver $Q_{f}.$ Assume  $Q_{f}$ is not the round trip
quiver. Let $(f,\delta)$ be an interlacing of $k^{E}$ determined
by a coloration $a$ of $Q_{f}.$ The twisted tensor algebra
$k^{E}\otimes_{(f,\delta)} \left[k[X]/(X^{2} - X)\right]$ is
isomorphic to the radical square zero algebra $kQ_{(f,\delta)} /
\langle Q_{(f,\delta)}^{2}\rangle.$
\end{thm}

\begin{rem}
Both algebras have the same dimension since  $Q_{(f,\delta)}$ has
more vertices but less arrows  in the same quantity, namely the
loops.
\end{rem}

 {\bf Proof:} Note first that an algebra morphism
$$\varphi : k Q_{(f,\delta)}\ \longrightarrow\ k^{E}\otimes_{(f,\delta)}\left[k[X]/(X^{2}
- X)\right]$$ is determined by a coherent choice of images of the
vertices and the arrows of $Q_{(f,\delta)}$ since $k
Q_{(f,\delta)}$ is a tensor algebra on the vertices set algebra of
the arrows bimodule.

We define $\varphi (e) = e$ for non loop vertices of $Q_{f}$. In
case of a loop vertex $\ell$ of $Q_f$ we put
$$\varphi(\ell_{0}) = \ell X \mbox{ and } \varphi(\ell_{-1}) =
\ell(1 - X).$$

 We have to
verify that the set $\left(E\setminus\!\{\ell\}\right)\bigsqcup \{
\ell X, \ell(1-X)\}$ is a complete set of orthogonal idempotents.
Recall that $k^{E} \otimes_{(f,\delta)} k[X]/(X^{2}-X)$ is
identified with $k^{E}[X]/(X^{2}-X)$ where the product is given by
the twist $Xb = \delta(b) + f(b)X.$ More precisely if $e \in E$,
we have
  $$X e = a(f(e)-e) + f(e)X
   = \sum_{\varphi(x) = e} a_{x} x - a_{e} e + \left(\sum_{\varphi(x) = e}
   x\right) X$$
If $e$ is $-1$-polarised, then
$$X e = e + \Bigl(\sum_{\varphi(x) = e}x\Bigr) X.$$
If $e$ is $0$-polarised
$$X e = \left(\sum_{\varphi(x) = e}x\right)\left(- 1 +X\right).$$
Finally if $e = \ell$ is a loop vertex, then
$$X \ell = -
\sum_{{\frac{}{}{\varphi(x) = \ell}\ {a_{x} = - 1}}} x
  +
  \left(
  \sum_{{\varphi(x) = \ell}} x
  \right) X$$

Note that

  $$(\ell X)^2 =  \ell X\ell X
   = \ell\left(-
\sum_{\varphi(x) = \ell \  a_{x} = - 1} x
  +
  \left(
  \sum_{{\varphi(x) = \ell}} x
  \right) X\right)X\\
  =\ell X$$

since $\varphi(x)=\ell$ and $\ell$ belongs to the sum
$\sum_{\varphi(x)=\ell}x$.

Similarly $\left[\ell(1-X)\right]^2=\ell(1-X).$ We also have that
$\left[\ell X\right]\left[\ell(1-X)\right]=0=
\left[\ell(1-X)\right]\left[\ell X\right]$ since $\ell
X\ell(1-X)=\ell X(1-X) =\ell X - \ell X = 0$.

Finally note that for a vertex $e\neq \ell$ we have $e\ell X=0$
and $\ell Xe=0$ for both possible colorations of $e$ since
$$1=\sum_{e\in E\setminus\{\ell\}} e + \ell X + \ell(1-X).$$

We consider now arrows of $Q_{(f,\delta)}$. For $e$ a non loop
vertex and $x\in E$ such that $\varphi(x)=e$ let ${}_xa_e$ be the
corresponding arrow. We assert that
$\varphi\left({}_xa_e\right)=xXe$ is a coherent choice since
$\varphi\left({}_xa_e\right)= \varphi\left(x\ {}_xa_e\ e\right)=
\varphi(x)\varphi\left({}_xa_e\right)\varphi(e).$ For $\varepsilon
= 0 \mbox{ or } -1$ let $x$ be a $\varepsilon$-polarized vertex of
$E$ such that $\varphi(x)=\ell$ and let ${}_xa_{\ell_\varepsilon}$
be the corresponding arrow. We put
$$\varphi\left({}_xa_{\ell_\varepsilon}\right)=xX\ell(1+\varepsilon X).$$
In order to prove that $\varphi$ is surjective note that for a non
loop vertex $e$ we have $yXe=0$ if $\varphi(y)\neq e$. Moreover,
if $\ell$ is a loop vertex and $y$ is not $\varepsilon$ polarized
then $yX\ell(1+\varepsilon X)=0$.

Since $$X\ =\ 1X1\ =\sum_{u,v \mbox{ \tiny  vertices of
}Q_{(f,\delta)}}vXu$$ we have that
$$\varphi\left(\sum_{\mbox{\tiny arrows
of }Q_{(f,\delta)}}a\right)=X.$$ Of course the entire algebra
$k^E$ is also in the image of $\varphi$, hence $\varphi$ is
surjective. An easy computation shows that a path of length two of
$Q_{(f,\delta)}$ has zero image. Since both algebras have the same
dimension the surjective map induces an isomorphism.

\begin{thm}

Let $E=\{u,v\}$ and let $f : k^{E}\longrightarrow k^{E}$ be an
algebra endomorphism given by a set-map $\varphi$ with quiver
$Q_{f} = u\cdot\!\leftrightarrows\!\cdot v$. Let $(f,\delta)$ be
an interlacing of $k^{E}$ determined by a coloration $a=a_uu+a_vv$
of $Q_{f}$, that is $a_u+a_v+1=0.$ The twisted tensor algebra
$k^{E}\otimes_{(f,\delta)} k[X]/(X^{2} - X)$ is isomorphic to the
radical square zero algebra $kQ_f / \langle Q_f^{2}\rangle.$

\end{thm}

\begin{rem}
The result of the Theorem is independent of the coloration, in
other words different colorations of the round trip quiver
provides the same twisted tensor product.
\end{rem}

 {\bf Proof:} The theory we have developed shows that the twisted
tensor product is the algebra $k^E[X]/(X-X^2)$ where the
interlacing is given by $Xu=a_vv-a_uu+vX$ and $Xv=a_uu-a_vv+uX$.
Recall that $a_u+a_v+1=0$. Note also that a direct computation can
be performed to show that the condition is necessary and
sufficient in order to have a well defined associative product.

We assert that the vertices of the algebra are $u$ and $v$ while
the arrows are $vXu$ and $uXv$. In other words we construct a
surjective  algebra map $\psi : kQ_f\longrightarrow
k^E[X]/(X-X^2).$ In order to prove this we have to show that $X$
is in the image of $\psi$. Indeed $X= uXu+vXv+vXu+uXv$, the last
two terms are in the image of $\psi$ by construction while
$uXu=-a_uu$ and $vXv=-a_vv$. Finally note that $vXuXv=0=uXvXu$
which shows that $\psi$ provides a surjective map between algebras
of dimension $4$, consequently they are isomorphic.

It is well known and easy to prove that a complet invariant up to
isomorphism of radical square zero algebras $kQ/<Q^2>$ is provided
by the quiver $Q$. Of course this fact will be a crucial piece in
order to classify the non-commutative duplicates of a finite set
through the related quiver.

We already note that a connected one-valued quiver not reduced to
a loop, which is not the round trip quiver and with even length
oriented cycle has exactly two colorations. In case there is no
loop the related quiver do not depend on the coloration. In
contrast in case there is a loop the related quivers differs
depending on the coloration.

As we quote before, the related quiver of a connected quiver
containing a loop has two connected components. Both have
precisely one sink vertex and all other vertices are one-valued.
We call such quivers \emph{one-sink-one-valued}, observe that a
quiver reduced to a vertex belongs to this family.

This discussion proves the following

\begin{thm}
Let $E$ be a finite set. The complete list up to isomorphism of
non-commutative duplicates of $E$ are provided by the disjoint
union of quivers as follows
\begin{itemize}
\item connected one-valued quivers with an even length proper
oriented cycle, \item connected one-sink-one-valued quivers
\end{itemize}
such that the total number of vertices involved equals $\mid
\!E\!\mid\!+ \ L_1$ where $L_1$ is the number of connected
one-sink-one-valued components (corresponding to the number of
loops of the original quiver).
  The non-commutative duplicates are the
radical square zero algebras obtained with the opposites of those
quivers.
\end{thm}

Note that different $2$-interlacings can provide the same related
quiver. Coloration switches on non loop connected component
provides the same related quiver. In case the set map has no fixed
points (that is the quiver $Q_f$ has no loops), the fiber of the
related quiver is precisely described this way. Its cardinality is
$2$ raised to the number of such connected components.

\section{Hochschild (co)homology and cyclic homology}

In order to compute the Hochschild cohomology of non-commutative
duplicates of a finite set, we will use the classification of the
preceding section -- we have shown that they are provided by a
precise family of quivers and the corresponding radical square
zero algebras. We will also use the main results from \cite{ci}
where the dimensions of the Hochschild cohomology for these
algebras are computed.

Let $Q$ be a connected quiver provided by a finite set of vertices
$Q_0$, a finite set of arrows $Q_1$ and two maps $s,t:
Q_1\rightarrow Q_0$ determining the source and the target of each
arrow. Let $Q_n$ be the set of paths of length $n$ in $Q$, namely
the sequences of $n$ concatenated arrows. We denote $(kQ)_2$ the
quotient of the path algebra of $Q$ modulo the two sided ideal
generated by $Q_2$.

The dimension of the Hochschild cohomology of $(kQ)_2$ involves
the sets of \emph{parallel paths}. More precisely two paths are
said to be parallel if they share the same source and the same
terminal vertices. In case $X$ and $Y$ are sets of paths,
$X/\!\!/Y$ is the set of couples of parallel paths
$(\alpha,\beta)$ where $\alpha\in X$ and $\beta\in Y$.

The following result is proved in \cite{ci}. It concerns any
connected quiver different from a crown, where a $c$-crown is the
quiver with $c$ vertices cyclically labelled by the cyclic group
of order $c$, and $c$ arrows $a_0,\dots , a_{c-1}$ such that
$s(a_i) = i$ and $t(a_i) = i + 1$.

\begin{thm}

Let $Q$ be a connected quiver which is not a crown. Then

$$\dim HH^n((kQ)_2)=\  \mid\! (Q_n/\!\!/Q_1)\mid\ -\ \mid
(Q_{n-1}/\!\!/Q_0)\mid  \mbox{ in case } n\geq 2$$

$$\dim HH^1((kQ)_2)=\ \mid\! (Q_1/\!\!/Q_1)\mid\ - \ \mid Q_0 \mid\ + 1$$

$$\dim HH^0((kQ)_2)=\ \mid\! (Q_1/\!\!/Q_0)\mid\ +
1$$
\end{thm}

\begin{rem}
In \cite{ci} there is a one unity error in the computation of the
dimension of the degree one Hochshchild cohomology.
\end{rem}

The preceding Theorem has an important consequence, see Corollary
3.2 of \cite{ci}:

\begin{cor}
Let $Q$ be a connected quiver which is not a crown. The graded
cohomology $HH^*((kQ_2))$ is finite dimensional if and only if $Q$
has no oriented cycles.
\end{cor}

Observe that in case of a one-valued quiver the set $\mid\!
(Q_n/\!\!/Q_1)\mid$ for $n\geq 2$ is only originated by oriented
cycles. More precisely an oriented cycle of length $n-1$ preceded
by an arrow provides the set of couples in $\mid\!
(Q_n/\!\!/Q_1)\mid $ for $ n\geq 2 $. Observe also that in case of
a connected one-sink-one-valued quiver $\mid\! (Q_n/\!\!/Q_1)\mid=
\emptyset$ for $n\geq 2$.

\begin{thm}
Let $A$ be a non-commutative duplicate algebra of a finite set,
given by
 a disjoint union of connected one-sink-one-valued quivers only.
Equivalently, there is no oriented cycles other than loops in the
quiver of the set map defining the non-commutative duplicate.

Then $HH^n(A)=0$ for $n\geq 2$, while $\dim HH^0(A)$ is the number
of connected components of the quiver and $\dim HH^1(A)$ is the
Euler characteristic of the underlying graph of $Q$.

\end{thm}

\begin{thm}
Let $A$ be a non-commutative duplicate algebra of a finite set
provided by a $2$-nilpotent algebra associated to a quiver which
has an oriented cycle. Then the graded cohomology $HH^*(A)$ is
infinite dimensional.
\end{thm}
 {\bf Proof:} The result is clear from the previous considerations,
and from the fact that a crown has infinite dimensional
cohomology, see Proposition 3.3 of \cite{ci}.

Finally we resume our results in order to make a direct link with
the set $Y_{k^E}$ of $2$-interlacings.

\begin{thm}
Let $(f,\delta)$ be a $2$-interlacing of a finite set $E$, namely
$f$ is an endomorphism of $k^E$ provided by a set map $\varphi :
E\rightarrow E$ and $\delta : k^E\rightarrow {}^f\!(k^E)$ is an
idempotent derivation verifying $f = f^2 + \delta f + f\delta$.
Let $A$ be the corresponding non-commutative duplicate of $E$
obtained through the twisting tensor product, $A= k^E
\otimes_{(f,\delta)} \left[k[X]/(X^2-X)\right]$.

Then $H^*(A)$ is finite dimensional if and only if the quiver of
the set map $\varphi$ do not contains oriented cycles other than
possible loops. In that case the Hochschild cohomology vanishes in
degrees larger than $2$.

\end{thm}

 {\bf Proof:} The related quiver described in the previous section
eliminates loops of $Q_f$ without creating oriented cycles,
regardless of the coloration originated by the derivation.

Concerning cyclic homology we will use the results obtained in
\cite{ci2}. Recall that the \emph{circuit} associated to an
oriented cycle $\gamma$  is the orbit of $\gamma$  through the
natural action of the cyclic group of the same order than the
length of $\gamma$. The set of circuits of length $j$ is denoted
$\Theta_j$.

Note that in case of a proper cycle the action of the
corresponding cyclic group is free on its orbit. We denote
$\Omega_a$ the set of proper circuits of length $a$. The next
result is proven in \cite[p. 139]{ci2}:

\begin{prop}
Let $Q$ be a quiver and $k$ be a field of characteristic zero.

For $n$ even
$$\dim HC_n\left((kQ_2)\right)= \ \mid\!\Theta_{n+1}\mid +\mid Q_0\mid.$$

For $n$ odd
$$\dim HC_n\left((kQ_2)\right)= \sum_{2\mid a \mid
(n+1)}\mid\Omega_a\mid.$$
\end{prop}

\begin{thm}

Let $k$ be a field of characteristic zero, let $E$ be a finite
set, let $A$ be a non-commutative duplicate algebra of $E$ and let
$\varphi: E\rightarrow E$ be the set algebra map corresponding to
the interlacing. The dimension of the even degree cyclic homology
of $A$ is constant regardless of the derivation:
$$\dim HC_{\mbox{\small \rm
even}}(A)=\  \mid\! E\mid + \mid \mbox{\rm loops}\mid.$$
\end{thm}

 {\bf Proof:} We know that there are no odd length oriented cycles in
the quiver of $f$, consequently the dimension is reduced to the
number of vertices of the related quiver. Recall that each loop in
the quiver of $f$ produces a new vertex in the related quiver.

In order to compute the odd cyclic homology, consider a
$2$-interlacing $(f,\delta)$ and the quiver $Q_f$ of $f$.  Let
$h(a)$ be the number of connected components  of $Q_f$ containing
a proper cycle of length $a$. We know that if $a$ is odd then
$h(a)=0$.

\begin{thm}
Let $k$ be a field of characteristic zero, let $E$ be a finite
set, let $A$ be a non-commutative duplicate algebra of $E$ and let
$f: E\rightarrow E$ be the set algebra map corresponding to the
$2$-interlacing. Then for $n$ odd
$$\dim HC_{n}(A)=\sum_{a\mid (n+1)}h(a).$$
\end{thm}

Finally the Hochschild homology can also be obtained using the
computation in \cite[p. 140]{ci2}.

\begin{thm}
Let $k$ be a field of characteristic different from $2$, let $E$
be a finite set, let $A$ be a non-commutative duplicate algebra of
$E$ and let $f: E\rightarrow E$ be the set algebra map
corresponding to the $2$-interlacing.

For $n$ odd

$$\dim HH_n (A) = \sum_{a\mid(n+1)} h(a).$$

For $n$ even and $n\neq 0$

$$\dim HH_n (A) = \sum_{a\mid n} h(a).$$

For $n=0$

$$\dim HH_0(A) = \  \mid\! E\mid + \mid \mbox{\rm loops}\mid.$$
\end{thm}



\end{document}